\documentclass[11pt,a4paper]{amsart}
\usepackage[english]{babel}
\usepackage[T1]{fontenc}
\usepackage[latin1]{inputenc}
\usepackage{amsfonts}
\usepackage{amssymb}
\usepackage{mathrsfs}
\usepackage{amsmath}
\usepackage{amsthm}
\usepackage{latexsym}
\usepackage{graphicx}
\usepackage{pst-all}
\newtheorem{Lemm}{Lemma}[section]
\newtheorem{Prop}{Proposition}[section]
\newtheorem{Corol}{Corollary}[section]
\newtheorem{rema}{Remark}[section]
\newtheorem{defi}{Definition}[section]

\def\qed{\hfill$\square$\par \bigskip}
\newenvironment{Demo}[1]{{\bf Proof #1~}}{\qed}
\newcommand{\R}{\mathbb{R}}

\newcommand{\N}{\mathbb{N}}

\DeclareMathOperator{\Ima}{Im}
\newcommand{\norm}[1]{\left\Vert#1\right\Vert}
\newcommand{\abs}[1]{\left\vert#1\right\vert}

\usepackage{graphicx}

\newcommand{\MSC}{MSC{}}
\newcommand{\Keywords}{Keywords{}}

\usepackage{times}

\begin{document}

\author{Hèla Ayadi}
\title{Semi-Fredholmness of the discrete Gauss-Bonnet Operator}
\address{ Unité de recherche Mathématiques et applications (UR/13ES47) à la faculté des sciences de Bizerte. Laboratoire De Math\'{e}matiques Jean Leray, Universit\'{e} de Nantes.}
\email{halaayadi@yahoo.fr. Hela.Ayadi@univ-nantes.fr.}

\begin{abstract}
In the context of an infinite locally finite weighted graph, we give a necessary and sufficient condition for semi-Fredholmness of the Gauss-Bonnet operator. This result is a discrete version of the theorem of Gilles Carron in the continuous case \cite{CARR}. In addition, using a criterion of Anghel \cite{AN}, we give a sufficient condition to have an operator of Gauss-Bonnet with closed range. Finally, this work can be considered as an extension of the work of Colette Ann\'{e} and Nabila Torki-Hamza \cite{ytt}.
\end{abstract}
\date{\today,\\
\MSC.\;{39A12, 05C63, 05C22, 47A53.} \\
\Keywords.\;{Infinite weighted graph, Discrete Gauss-Bonnet Operator, 
Non-parabolicity at Infinity, Closed Range Operator.}\\
}
\maketitle
\tableofcontents

\section{Introduction}
 Dirac type operators have become of central importance in many branches of mathematics such as PDE's, differential geometry and topology (see \cite{dirac}, \cite{dirac1}, \cite{dirac2}..), since the introduction in 1928 by the physicist Paul Dirac of a first-order linear differential operator whose square is the Laplacian operator. In particular, this paper focuses on the conditions to have semi-Fredholmness of the discrete Gauss-Bonnet operator needed to approach the Hodge decomposition theorem \cite{ytt}. In fact, we present a discrete version of the work of G. Carron \cite{CARR}, which defines a new concept "non-parabolicity at infinity" to have the Gauss-Bonnet operator with closed range. Indeed, G. Carron's condition  is quite weaker than the one given by Anghel \cite{AN}. Moreover, we provide a new sufficient condition to obtain a Gauss-Bonnet operator semi-Fredholm. Finally, we give two explicit examples one example verifying the property of non-parabolicity at infinity, and the other not.

\section{Preliminaries}
\subsection{Definitions and notations}

\begin{itemize}

 \item A graph $G$ is a couple $(\mathcal{V}, \mathcal{E})$ where $\mathcal{V}$ is a set at most countable whose elements are called vertices and $\mathcal{E}$ is a set of oriented edges, considered as a subset of $\mathcal{V}\times\mathcal{V}$.
 \item  If the graph $G$ has a finite set of vertices, it's called a finite graph. Otherwise, $G$ is called infinite graph.

 \item We assume that $\mathcal{E}$ is symmetric without loops:
$$v \in \mathcal{V} \Rightarrow (v, v) \notin \mathcal{E},
 \;\; (v_{1}, v_{2})\in \mathcal{E} \Rightarrow (v_{2}, v_{1})\in \mathcal{E}.$$

 \item  Choosing an orientation of $G$ consists of defining a partition of $\mathcal{E}$: $\mathcal{E}^{+}\sqcup \mathcal{E}^{-}=\mathcal{E}$
 $$(v_{1}, v_{2})\in \mathcal{E}^{+} \Leftrightarrow (v_{2},v_{1})\in \mathcal{E}^{-}.$$

 \item For $e=(v_{1}, v_{2})$, we denote
$$e^{-}=v_{1},\; e^{+}=v_{2}\; and\; -e=(v_{2}, v_{1}).$$

 \item The graph $G$ is connected if, any two vertices $x$, $y$ in $\mathcal{V}$ can be joined by a path of edges $\gamma_{xy}$, that means,
     $$\gamma_{xy}=\{e_{k}\}_{k=1,..., n}\; with\; e^{-}_{1}=x,\; e^{+}_{n}=y\; and\; if\; n\geq2\; ,\;\forall j \; ;1\leq j \leq (n-1) \Rightarrow e^{+}_{j}=e^{-}_{j+1}.$$

 \item  The degree (or valence) of a vertex $x$ is the number of edges emanating from $x$. We denote
$$ deg(x):= \sharp\{  e\in \mathcal{E};\;e^{-} = x  \}.$$

 \item If $deg(x)< \infty,\;\forall x \in \mathcal{V}$,  we say that $G$ is a locally finite graph.

 \end{itemize}

 \subsection{The weighted graph}

The weighted graph $(G, c, r)$ is given by the graph
 $G=(\mathcal{V}, \mathcal{E})$, a weight on the vertices
 $c:\mathcal{V} \rightarrow ]0, \infty[$ and a weight on the edges
  $r: \mathcal{E} \rightarrow ]0, \infty[$such that $r(-e)=r(e)$.\\

\textbf{Examples:}
- An infinite electrical network is a weighted graph $(G,c , r)$ where the weights of the edges called resistances $r$; their reciprocals are called conductances. And the weights of the vertices given by $c(x)=\sum_{y\in \mathcal{V}} \frac{1}{ r(x, y)}< \infty,\; \forall x \in \mathcal{V}$.\\

-The graph $G$ called a simple graph where the weights of the edges and the vertices equals $1$.\\

\textbf{All the graphs we shall consider on the sequel will be weighted, connected and locally finite.}

\subsection{The notion of subgraph}

A subgraph of a graph $G$ is a graph $G_{K}:=(K, \mathcal{E}_{K})$ such that $K\subset\mathcal{V}$ and $\mathcal{E}_{K}:=\{e\in \mathcal{E};\; e^{-},\;e^{+} \in K\}.$

For such a subgraph we define:

\begin{itemize}
\item the vertex boundary : $$\partial K:=\{ x \in \mathcal{V} \setminus{K}; \exists y \in K, \; (x, y) \in   \mathcal{E}                      \},$$
\item the edge boundary: $$\partial \mathcal{E}_{K}:=\{ e \in \mathcal{E};\; e^{-} \in K\; and \;e^{+} \notin K \; or \; e^{+} \in K\; and \;e^{-} \notin K  \}.$$
\end{itemize}

\subsection{Functional spaces}
We denote the set of real functions on $\mathcal{V}$ by:
$$\mathcal{C}(\mathcal{V})=\{f:\mathcal{V} \rightarrow \R\}$$

and the set of functions of finite support by $\mathcal{C}_{0}(\mathcal{V})$.\\

Moreover, we denote the set of real skewsymmetric functions on $\mathcal{E}$ by:
$$\mathcal{C}^{a}(\mathcal{E})=\{\varphi:\mathcal{E} \rightarrow \R \;;
 \varphi(-e)=- \varphi(e)\}$$

and the set of functions of finite support by $\mathcal{C}^{a}_{0}(\mathcal{E})$.\\

We define on the weighted graph $(G, c, r)$ the following function spaces endowed of the scalar products.

\begin{description}
  \item[a)] $$l^{2}(\mathcal{V}):=\left\{ f\in \mathcal{C}(\mathcal{V})  ;\;
   \sum_{x\in \mathcal{V}} c(x) f^{2}(x)<\infty \right\},$$

with the inner product
$$\left<f, g\right>_{\mathcal{V}}=\sum_{x\in \mathcal{V}} c(x) f(x)g(x)$$

and the norm
$$\norm{f}_{l^{2}(\mathcal{V})}=\sqrt{\left<f, f\right>_{\mathcal{V}}}.$$

\item[b)] $$l^{2}(\mathcal{E}):=\left\{ \varphi \in \mathcal{C}^{a}(\mathcal{E});\;
\frac{1}{2}\sum_{e\in \mathcal{E}} r(e) \varphi^{2}(e)<\infty \right\},$$

with the inner product
$$\left<\varphi, \psi\right>_{\mathcal{E}}=
\frac{1}{2} \sum_{e\in \mathcal{E}} r(e) \varphi(e)\psi (e)$$

and the norm
$$\norm{\varphi}_{l^{2}(\mathcal{E})}=
\sqrt{\left<\varphi, \varphi\right>_{\mathcal{E}}}.$$

\end{description}

As a consequence, we define the direct sums of $l^{2}(\mathcal{V})$ and $l^{2}(\mathcal{E})$ by:
$$l^{2}(G):=l^{2}(\mathcal{V}) \oplus l^{2}(\mathcal{E})=\left\{(f, \varphi),\;f\in l^{2}(\mathcal{V}) \hbox{ and } \varphi \in l^{2}(\mathcal{E})\right\},$$
with the norm
$$\norm{(f, \varphi)}_{l^{2}(G)}^{2}:=
\norm{f}_{l^{2}(\mathcal{V})}^{2}+\norm{\varphi}_{l^{2}(\mathcal{E})}^{2}.$$
\subsection{Operators and properties}

\underline{The difference operator:} it is the operator
$$\mathrm{d}:\mathcal{C}_{0}(\mathcal{V})\longrightarrow \mathcal{C}^{a}_{0}(\mathcal{E}),$$

given by
$$\mathrm{d}(f)(e)=f(e^{+})-f(e^{-}).$$

\underline{The coboundary operator:}
it is $\delta$ the formal adjoint of $\mathrm{d}$. Thus it satisfies
\begin{equation}\label{op div}
\left<\mathrm{d}f, \varphi \right>_{\mathcal{E}}=
\left<f, \delta \varphi \right>_{\mathcal{V}}\end{equation}
for all $ f \in \mathcal{C}_{0}(\mathcal{V})$ and for all $ \varphi \in \mathcal{C}^{a}_{0}(\mathcal{E}).$\\

As consequence, we have the following formula characterizing $\delta$ :

\begin{Lemm}\label{lem1}
The coboundary operator $\delta$ is characterized by the formula

$$\delta \varphi (x)=\frac{1}{c(x)}  \sum_{e, e^{+}=x} r(e) \varphi (e),$$

for all $ \varphi \in \mathcal{C}^{a}_{0}(\mathcal{E}).$

\end{Lemm}

\begin{Demo}

For $f \in \mathcal{C}_{0}(\mathcal{V})$ and $\varphi \in \mathcal{C}_{0}^{a}(\mathcal{E})$,
using (\ref{op div}), we get
\begin{eqnarray*}
\left<\mathrm{d}f, \varphi \right>_{\mathcal{E}} &=& \frac{1}{2} \sum_{e\in \mathcal{E}} r(e) \mathrm{d}f(e) \varphi(e)\\
&=& \frac{1}{2} \sum_{e\in \mathcal{E}} r(e) \left(f(e^{+})-f(e^{-}) \right)\varphi(e)\\
&=&\frac{1}{2} \sum_{x\in \mathcal{V}} f(x) \left(  \sum_{e, e^{+}=x}     r(e)\varphi(e)-\sum_{e, e^{-}=x}     r(e)\varphi(e)\right).
\end{eqnarray*}

But, $r(-e)=r(e)$ and  $\displaystyle{\sum_{e, e^{+}=x}     r(e)\varphi(e)=-\sum_{e, e^{-}=x}     r(e)\varphi(e)}$.\\

So we have,
\begin{eqnarray*}
\left<\mathrm{d} f, \varphi \right>_{\mathcal{E}}&=&\sum_{x\in \mathcal{V}} c(x) f(x) \left(\frac{1}{c(x)} \sum_{e, e^{+}=x}     r(e)\varphi(e)\right)\\
&=&\left<f, \delta \varphi \right>_{\mathcal{V}}.
\end{eqnarray*}

\end{Demo}

We introduce now a very important result inspired by \cite{Nyam}.
\begin{Lemm}\label{lem2}
Let $x$ and $ x_{0}$ in $\mathcal{V}$, then there exists a positive constant $C_{xx_{0}}$ such that
\begin{equation}\label{eqlem2}
\abs{f(x)}\leq C_{xx_{0}} \left( \abs{f(x_{0})} + \norm{\mathrm{d}f}_{l^{2}(\mathcal{E})}\right),
\end{equation}
for all $ f \in \mathcal{C}_{0}(\mathcal{V}).$
\end{Lemm}

\begin{Demo}

As $G$ is connected, then we can find a path $\gamma_{xx_{0}}$ joining $x$ to $x_{0}$, i.e,
$$\gamma_{xx_{0}}=\{e_{k}\}_{k=1,..., n}\; with\; e^{-}_{1}=x,\; e^{+}_{n}=x_{0}\hbox{ and if } n\geq2\; ,\;\forall j \; ;1\leq j \leq (n-1) \Rightarrow e^{+}_{j}=e^{-}_{j+1}.$$

Then, using the triangle inequality, we have
\begin{eqnarray*}
\abs{f(x)-f(x_{0})}&=&\abs{f(x)-f(e^{+}_{1})+f(e^{+}_{1})-f(e^{+}_{2})+...+f(e^{+}_{n-1})-f(x_{0})}\\
&\leq&\abs{\mathrm{d}f(e_{1})}+\abs{\mathrm{d}f(e_{2})}+...+\abs{\mathrm{d}f(e_{n})}\\
&\leq& \sum_{e\in \gamma_{xx_{0}}} \frac{1}{\sqrt{r(e)}} \sqrt{r(e)}\abs{\mathrm{d}f(e)}.
\end{eqnarray*}
Applying the Cauchy-Schwarz inequality, we obtain
\begin{eqnarray*}
\abs{f(x)-f(x_{0})}&\leq& \left(\sum_{e\in \gamma_{xx_{0}}} \frac{1}{{r(e)}}\right)^{\frac{1}{2}} \left(\sum_{e\in \gamma_{xx_{0}}} {r(e)}(\mathrm{d}f(e))^{2}\right)^{\frac{1}{2}}\\
&\leq& S_{xx_{0}}\left(\sum_{e\in \mathcal{E}} {r(e)}(\mathrm{d}f(e))^{2}\right)^{\frac{1}{2}}\\
&\leq& S_{xx_{0}}\norm{\mathrm{d}f}_{l^{2}(\mathcal{E})},
\end{eqnarray*}
with $S_{xx_{0}}=\left(\displaystyle{\sum_{e\in \gamma_{xx_{0}}} \frac{1}{{r(e)}}}\right)^{\frac{1}{2}}.$\\

Thus, we deduce that
\begin{eqnarray*}
\abs{f(x)}&\leq& \abs{f(x)-f(x_{0})}+\abs{f(x_{0})}\\
&\leq& S_{xx_{0}}\norm{\mathrm{d}f}_{l^{2}(\mathcal{E})}+\abs{f(x_{0})}\\
&\leq& C_{xx_{0}}\left(\norm{\mathrm{d}f}_{l^{2}(\mathcal{E})}+\abs{f(x_{0})}\right),
\end{eqnarray*}
with $C_{xx_{0}}=max(S_{xx_{0}}, 1)$.
\end{Demo}

Before giving another important result, for $f \in \mathcal{C}_{0}(\mathcal{V})$, we define the mean value $\overline{f}$ of $f$ by
$$\overline{f}(e)=\frac{f(e^{+})+f(e^{-})}{2} $$

for all $e\in \mathcal{E}.$\\

And we have from \cite{mus} the following derivation property:

\begin{Lemm}\label{lem3}
For $f, g \in \mathcal{C}_{0}(\mathcal{V})$ and $\varphi \in \mathcal{C}^{a}_{0}(\mathcal{E})$, it follows
\begin{equation}\label{eqlem3}\mathrm{d}(fg)(e)=f(e^{+})\mathrm{d}g(e)+g(e^{-})\mathrm{d}(f)(e).\end{equation}
\begin{equation}\label{eqqlem3}\delta(\overline{f} \varphi)(x)= f(x) \delta \varphi(x)-\frac{1}{2 c(x)} \sum_{e, e^{+}=x} r(e)\mathrm{d}(f)(e) \varphi (e).\end{equation}
\end{Lemm}

\begin{Demo}

For $f, g \in \mathcal{C}_{0}(\mathcal{V})$ and $e\in \mathcal{E}$,
\begin{eqnarray*}
\mathrm{d}(fg)(e)&=& (fg)(e^{+})- (fg)(e^{-})\\
&=& f(e^{+})\left(g(e^{+})-g(e^{-})\right)+g(e^{-})\left(f(e^{+})-f(e^{-})\right)\\
&=& f(e^{+})\mathrm{d}(g)(e)+g(e^{-})\mathrm{d}(f)(e).
\end{eqnarray*}
On the other hand, for $\varphi \in \mathcal{C}^{a}_{0}(\mathcal{E})$ applying the characterization of $\delta$ from Lemma (\ref{lem1}) to the function $\overline{f} \varphi \in \mathcal{C}^{a}_{0}(\mathcal{E})$, we have
\begin{eqnarray*}
\delta(\overline{f} \varphi)(x)&=&\frac{1}{c(x)}  \sum_{e, e^{+}=x} r(e) (\overline{f}\varphi) (e)\\
&=&\frac{1}{c(x)}  \sum_{e, e^{+}=x} r(e) \left(\frac{f(e^{+})+f(e^{-})}{2}  \right) \varphi (e)\\
&=&\frac{1}{c(x)}  \sum_{e, e^{+}=x} r(e) f(e^{+})\varphi (e)+\frac{1}{c(x)}  \sum_{e, e^{+}=x} r(e) \left(\frac{f(e^{-})-f(e^{+})}{2}  \right) \varphi (e)\\
&=&f(x)\frac{1}{c(x)}  \sum_{e, e^{+}=x} r(e) \varphi (e) +\frac{1}{2c(x)} \sum_{e, e^{+}=x} r(e) \mathrm{d}(f)(-e) \varphi(e)\\
&=&f(x)\delta(\varphi)(x) -\frac{1}{2c(x)} \sum_{e, e^{+}=x} r(e) \mathrm{d}(f)(e) \varphi(e).
\end{eqnarray*}
\end{Demo}

\underline{The Gauss-Bonnet operator:} it is the endomorphism

$$D=\mathrm{d}+\delta: \mathcal{C}_{0}(\mathcal{V}) \oplus \mathcal{C}^{a}_{0}(\mathcal{E})\longrightarrow \mathcal{C}_{0}(\mathcal{V}) \oplus \mathcal{C}^{a}_{0}(\mathcal{E})$$ with, $$D(f,\varphi)=\delta \varphi+\mathrm{d} f,\; \; \forall(f,\varphi)\in \mathcal{C}_{0}(\mathcal{V}) \oplus \mathcal{C}^{a}_{0}(\mathcal{E}).$$
And it is a symmetric operator.

\section{Non-parabolicity at infinity}
Now we introduce the discrete result of Carron \cite{CARR}:
\begin{defi}
 We say that $D$ is non-parabolic at infinity if there is a finite subgraph $G_{K}$ of $G$ such that for all finite subset $U$ of $G \setminus G_{K}$, there exists a positive constant $C=C({U})$ such that holds the following inequality
$$C \norm{(f,\varphi)}_{l^{2}(U)} \leq \norm{D(f,\varphi)}_{l^{2}(G\setminus G_{K})},\; \forall (f, \varphi) \in \mathcal{C}_{0}(\mathcal{V}\setminus K)\times \mathcal{C}^{a}_{0}(\mathcal{E}\setminus \mathcal{E}_{K}).$$

\end{defi}

\begin{rema}
We call a finite subset $U$ of $G$ a couple $U:=(\mathcal{V}_{U}, \mathcal{E}_{U})$ such that $\mathcal{V}_{U}$ is a finite subset of $\mathcal{V}$ and $\mathcal{E}_{U}$ is a finite subset of $\mathcal E$. And, we denote
$$\norm{(f, \varphi)}_{l^{2}(U)}^{2}=\norm{f}_{l^{2}(\mathcal{V}_{U})}^{2}+\norm{\varphi}_{l^{2}(\mathcal{E}_{U})}^{2}.$$
\end{rema}
\begin{defi}\label{neighb}

$G_{\widetilde{K}}$ is \textbf{a neighborhood} of $G_{K}$ if $G_{\widetilde{K}}:=({\widetilde{K}}, \mathcal{E}_{\widetilde{K}})$ is a finite subgraph of $G$ such that
$$
\left\{
  \begin{array}{ll}
  i)\; K \subset {\widetilde{K}}\; finite,\\
\\
  ii)\; \mathcal{E}_{K}\sqcup \partial \mathcal{E}_{K} \subset \mathcal{E}_{\widetilde{K}},\\
\\
  iii)\; e=(x, y) \in \mathcal{E}_{\widetilde{K}} \Rightarrow x,\; y \in {\widetilde{K}}.
   \end{array}
\right.
$$
Since we can define \textbf{the smallest neighborhood} of $G_{K}$ by  $G_{\widetilde{K}_{0}}$, where $G_{\widetilde{K}_{0}}$ is a finite subgraph of $G$ contains $G_{K}$ and its boundary.
\end{defi}
\begin{rema}
In \cite{keller}, $G_{\widetilde{K}_{0}}$ is called a combinatorial neighborhood of $G_{K}$.
\end{rema}
\begin{Lemm}\label{lem}
If $D$ is non-parabolic at infinity then, for every finite subset $U$ of $G$ there exists a positive constant $C'=C'(U)$ such that
\begin{equation}
C' \norm{(f, \varphi)}_{l^{2}(U)} \leq \norm{D(f, \varphi)}_{l^{2}(G)}+\norm{(f, \varphi)}_{l^{2}(G_{\widetilde{K}})},\;\forall (f, \varphi)\in \mathcal{C}_{0}(\mathcal{V}) \oplus \mathcal{C}^{a}_{0}(\mathcal{E}),
\end{equation}
where $G_{\widetilde{K}}$ is a neighborhood of $G_{K}$.
\end{Lemm}

\begin{Demo}

Since $U$ is a finite subset of $G$ it can be reduced to a point or an edge.\\

Let $x$ any vertex of $G$, we start by proving 
$$C' \abs{f(x)}\leq \norm{\mathrm{d}f}_{l^{2}(\mathcal{E})}+\norm{f }_{l^{2}(\widetilde{K})},\;\forall f\in \mathcal{C}_{0}(\mathcal{V}).$$
 
$G_{\widetilde{K}}$ is a finite subgraph of $G$, so according to Lemma \ref{lem2}, we obtain
\begin{equation}f^{2}(x) \leq C_{1} \left(\norm{f}_{l^{2}({\widetilde{K}})}^{2}+ \norm{\mathrm{d}f }_{l^{2}(\mathcal{E})}^{2}\right),\end{equation}
where $C_{1}$ is a positive constant which depends on $x$ and $\widetilde{K}$. Indeed:\\
let $x\in \mathcal{V}$ and $x_{0}\in \widetilde{K}$, using Lemma \ref{lem2}, we obtain
\begin{equation}f^{2}(x)\leq C_{xx_{0}} \left( {f^{2}(x_{0})} + \norm{\mathrm{d}f}^{2}_{l^{2}(\mathcal{E})}\right).\end{equation}
Multiplying (3.7) by $c(x_{0})>0$, we get
\begin{eqnarray*}
c(x_{0}){f^{2}(x)}&\leq& C_{xx_{0}} \left(c(x_{0}) {f^{2}(x_{0})} + c(x_{0})\norm{\mathrm{d}f}^{2}_{l^{2}(\mathcal{E})}\right)\\
&\leq& C_{xx_{0}} \left(\norm{f}_{l^{2}({\widetilde{K}})}^{2} + c(x_{0})\norm{\mathrm{d}f}^{2}_{l^{2}(\mathcal{E})}\right)\\
&\leq& C'_{xx_{0}} \left(\norm{f}_{l^{2}({\widetilde{K}})}^{2} + \norm{\mathrm{d}f}^{2}_{l^{2}(\mathcal{E})}\right),
\end{eqnarray*}
where $C'_{xx_{0}}=max(C_{xx_{0}}, c(x_{0})C_{xx_{0}})$.\\
Then, we have
 $$f^{2}(x)\leq \frac{C'_{xx_{0}}}{c(x_{0})} \left(\norm{f}_{l^{2}({\widetilde{K}})}^{2} + \norm{\mathrm{d}f}^{2}_{l^{2}(\mathcal{E})}\right).$$
Finally, we obtain
$$ f^{2}(x)\leq C_{1} \left(\norm{f}_{l^{2}({\widetilde{K}})}^{2} + \norm{\mathrm{d}f}^{2}_{l^{2}(\mathcal{E})}\right)$$
where $C_{1}=\frac{C'_{xx_{0}}}{c(x_{0})}$.\\

On the other hand, we want to show the following inequality, for any edge $e\in \mathcal{E}$
$$C'' \abs{\varphi(e)}\leq \norm{\delta\varphi}_{l^{2}(\mathcal{V})}+\norm{\varphi }_{l^{2}(\mathcal{E}_{\widetilde{K}})},\;\forall \varphi\in \mathcal{C}_{0}(\mathcal{E}).$$
For $e\in \mathcal{E}_{K} \subset \mathcal{E}_{\widetilde{K}}$ finite, we have

$$\varphi^{2}(e) \leq \norm{\varphi}_{l^{2}(\mathcal{E}_{\widetilde{K}})}^{2}\leq \norm{\varphi}_{l^{2}(\mathcal{E}_{\widetilde{K}})}^{2}+\norm{\delta\varphi}_{l^{2}(\mathcal{V})}^{2}.$$

And if $e\in \mathcal{E}\setminus \mathcal{E}_{K}$, we consider the indicator function of $K^{c}$, denoted by $\chi $
\begin{equation}\label{cut}\chi(x)=\left\{
  \begin{array}{ll}
  0 \;\hbox{if}\; x\in K\\
\\
    1\; \hbox{otherwise}.
   \end{array}
\right.
\end{equation}
which gives
\begin{equation*}d\chi(e)=\left\{
  \begin{array}{ll}
   0 \;\hbox{if}\; e\in \mathcal{E}_{K},\\
\\
   \pm 1\; \hbox{if} \;e\in \partial \mathcal{E}_{K},\\
   \\
   0\; \hbox{otherwise}.
   \end{array}
\right.
 \And \; \;\overline{\chi}(e)=\left\{
  \begin{array}{ll}
  0 \;\hbox{if}\; e\in \mathcal{E}_{K},\\
\\
   \frac{1}{2}\; \hbox{if} \;e\in \partial \mathcal{E}_{K},\\
   \\
   1\; \hbox{otherwise}.
   \end{array}\right.\end{equation*}

Let $\varphi \in \mathcal{C}^{a}_{0}(\mathcal{E})$, we have then $\overline{\chi}\varphi$ with finite support in $\mathcal{E}\setminus\mathcal{E}_{K}$. Thus, applying the definition of the non-parabolicity at infinity of $D$ to the function $(0, \overline{\chi}\varphi)$, we obtain
$$\norm{\overline{\chi}\varphi}_{l^{2}(U)}^{2}\leq C\norm{\delta (\overline{\chi}\varphi)}_{l^{2}(\mathcal{V})}^{2},$$
where $C=\frac{1}{C(U)}$.\\

Since we have $e \in \mathcal{E}\setminus{\mathcal{E}_{K}}$, this implies that
\begin{equation}\label{lema}\varphi^{2}(e)
\leq C  \norm{\delta (\overline{\chi}\varphi)}_{l^{2}(\mathcal{V})}^{2}.\end{equation}
The derivation property of Lemma (\ref{lem3}), gives
$$\delta(\overline{\chi} \varphi)(x)= \chi(x) \delta \varphi(x)-\frac{1}{2 c(x)} \sum_{e, e^{+}= x} r(e) \mathrm{d}(\chi)(e) \varphi(e).$$
And by the inequality $(a-b)^{2}\leq 2 (a^{2}+b^{2} )$, we obtain
\begin{eqnarray*}
\norm{\delta (\overline{\chi}\varphi)}_{l^{2}(\mathcal{V})}^{2}&=& \sum_{x\in \mathcal{V}} c(x) (\delta (\overline{\chi}\varphi))^{2}\\
&\leq& 2 \left[\underbrace{\sum_{x\in \mathcal{V}} c(x)\left(\chi(x)\delta \varphi(x)\right)^{2}}_{I}+\underbrace{\sum_{x\in \mathcal{V}} c(x)\left(\frac{1}{2 c(x)} \sum_{e, e^{+}= x} r(e) \mathrm{d}(\chi)(e) \varphi (e)\right)^{2}}_{J}\right].
\end{eqnarray*}
So, for the first term we have
\begin{equation}\label{sum1}I=\sum_{x\in \mathcal{V}\setminus K} c(x)\left(\delta \varphi(x)\right)^{2}\leq\norm{\delta\varphi}_{l^{2}(\mathcal{V})}^{2}\end{equation}
and for the second one, we get
\begin{equation}\label{sum2}
J=\underbrace{\sum_{x\in K} \frac{1}{2 c(x)}\left( \sum_{e, e^{+}= x} r(e) \mathrm{d}(\chi)(e) \varphi (e)\right)^{2}}_{J_{1}}+\underbrace{\sum_{x\in \mathcal{V}\setminus K} \frac{1}{2 c(x)}\left( \sum_{e, e^{+}= x} r(e) \mathrm{d}(\chi)(e) \varphi (e)\right)^{2}}_{J_{2}}.\end{equation}
Using that $\mathrm{supp} (\mathrm{d}\chi ) = \partial \mathcal{E}_{K}\subset \mathcal{E}_{\widetilde{K}}$ and  the Cauchy-Schwarz inequality, we obtain
\begin{eqnarray*}
J_1&=& \sum_{x\in K} \frac{1}{2 c(x)}\left( \sum \limits_{\underset{e\in \mathrm{supp} (\mathrm{d}\chi)}{e, e^{+}= x}} r(e) \varphi (e)\right)^{2}\\
&=& C_{K}\left( \sum_{e\in \mathrm{supp} (\mathrm{d}\chi)}r(e) \varphi (e)\right)^{2}\\
&\leq& C_{K} \left(\sum_{e\in \mathrm{supp} (\mathrm{d}\chi) } r(e)\right)\left(\sum_{e\in supp (d\chi) } r(e) \varphi ^{2}(e)\right)\\
&\leq& C_{K} C'_{K} \sum_{e\in \mathcal{E}_{\widetilde{K}} } r(e) \varphi ^{2}(e)\\
&=&C_{2}\norm{\varphi}_{l^{2}(\mathcal{E}_{\widetilde{K}})}^{2},
\end{eqnarray*}
where $C_{K}=\max \limits_{x \in K } \frac{1}{2 c(x)}$, $C'_{K}=\sharp \mathcal{E}_{\widetilde{K}}\max\limits_{e \in \mathcal{E}_{\widetilde{K}}} r(e)$ and $C_{2}=C_{K} C'_{K}$.\\

And for $J_2$, we have $e=(e^{-}, e^{+}) \in \mathrm{supp} (\mathrm{d}\chi)=\partial \mathcal{E}_{K}$, so if $ e^{-}\in K,\;e^{+} \in \partial K$.

\begin{eqnarray*}
J_2&=& \sum_{x\in \partial K} \frac{1}{2 c(x)}\left( \sum \limits_{\underset{e\in \mathrm{supp} (\mathrm{d}\chi)}{e, e^{+}= x}} r(e) \varphi (e)\right)^{2}\\
&=& C''_{K}\left( \sum_{e\in \mathrm{supp} (\mathrm{d}\chi)}r(e) \varphi (e)\right)^{2}\\
&\leq& C''_{K} \left(\sum_{e\in \mathrm{supp} (\mathrm{d}\chi) } r(e)\right)\left(\sum_{e\in \mathrm{supp} (\mathrm{d}\chi) } r(e) \varphi ^{2}(e)\right)\\
&\leq& C''_{K} C'_{K} \sum_{e\in \mathcal{E}_{\widetilde{K}} } r(e) \varphi ^{2}(e)\\
&=&C'_{2}\norm{\varphi}_{l^{2}(\mathcal{E}_{\widetilde{K}})}^{2},
\end{eqnarray*}

where $C''_{K}=\max \limits_{x \in \partial K } \frac{1}{2 c(x)}$ and $C'_{2}=C''_{K} C'_{K}$.\\

Thus, $(\ref{sum2})$ becomes
\begin{equation} \label{res sum2}
J\leq C''_{2}\norm{\varphi}_{l^{2}(\mathcal{E}_{\widetilde{K}})}^{2},
\end{equation}
where $ C''_{2}=max( C_{2}, C'_{2})$.\\

So by $(\ref{sum1})$ and $(\ref{res sum2})$, we get
\begin{equation}\label{lemb}
\norm{\delta (\overline{\chi}\varphi)}_{l^{2}(\mathcal{V})}^{2}\leq max(2 , 2 C''_{2})\left(\norm{\delta\varphi}_{l^{2}(\mathcal{V})}^{2}+ \norm{\varphi}_{l^{2}(\mathcal{E}_{\widetilde{K}})}^{2}\right).
\end{equation}
Finally, (\ref{lema}) and (\ref{lemb}) give
$$\varphi^{2}(e)\leq \widetilde{C}\left(\norm{\delta\varphi}_{l^{2}(\mathcal{V})}^{2}+ \norm{\varphi}_{l^{2}(\mathcal{E}_{\widetilde{K}})}^{2}\right)$$
where $\widetilde{C}=\frac{2 max(1, C''_{2})}{C}$.

\end{Demo}

\begin{Prop}
If $D$ is non-parabolic at infinity, then we can construct a Hilbert space $W$ such that :
\begin{enumerate}
  \item$\mathcal{C}_{0}(\mathcal{V}) \oplus \mathcal{C}^{a}_{0}(\mathcal{E})$ is dense in $W$.

\item The injection of $\mathcal{C}_{0}(\mathcal{V}) \oplus \mathcal{C}^{a}_{0}(\mathcal{E})$ to $\mathcal{C}(\mathcal{V}) \oplus \mathcal{C}^{a}(\mathcal{E})$  extends by continuity to $W$.

  \item $D:W \longrightarrow l^{2}(G)$ is a bounded operator.
  \end{enumerate}
  \end{Prop}
\begin{rema} 
In 1) and 2) we use the topology of ponctual convergence on $\mathcal{C}(\mathcal{V}) \oplus \mathcal{C}^{a}(\mathcal{E})$, it means, the sequence $(f_{n}, \varphi_{n})$ converges ponctually to $(f, \varphi)$ on $\mathcal{C}(\mathcal{V}) \oplus \mathcal{C}^{a}(\mathcal{E})$ if 
$f_{n}(x)$ converges to $f(x)$, $\forall x \in \mathcal{V}$ and $\varphi_{n}(e)$ converges to $\varphi(e)$, $\forall e \in \mathcal{E}$.
\end{rema} 
\begin{rema}
In Carron's paper \cite{CARR}, the injection of the space of functions with compact support to $l^{2}_{loc}$ extends by continuity to $W$. But, in our case we didn't need to introduce the space $l^{2}_{loc}$ because in discrete case this notion is trivial.
\end{rema}

 \begin{Demo}

Let us denote by $W$ the closure of $\mathcal{C}_{0}(\mathcal{V}) \oplus \mathcal{C}^{a}_{0}(\mathcal{E})$ for the norm

$$N_{\widetilde{K}}(f, \varphi)=\left(\norm{(f, \varphi)}_{l^{2}(G_{\widetilde{K}})}^{2}+\norm{D(f, \varphi)}_{l^{2}(G)}^{2}\right)^{\frac{1}{2}},$$

where $G_{\widetilde{K}}$ is a neighborhood of $G_{K}$ (see Definition $(\ref{neighb}$)).\\

\underline{\textbf{Aim i)}}: $N_{\widetilde{K}}$ is a norm on $W$, we just look at the nullity, we have
\begin{eqnarray*}
N_{\widetilde{K}}(f, \varphi)=0 &\Leftrightarrow & \norm{(f, \varphi)}_{l^{2}(G_{\widetilde{K}})}=0\hbox{ and }\norm{D(f, \varphi)}_{l^{2}(G)}=0\\
&\Leftrightarrow & \norm{f}_{l^{2}({\widetilde{K}})}=0, \;\norm{ \varphi}_{l^{2}(\mathcal{E}_{\widetilde{K}})}=0, \;\norm{\mathrm{d}f }_{l^{2}(\mathcal{E})}=0\hbox{ and }\norm{\delta \varphi }_{l^{2}(\mathcal{V})}=0.
\end{eqnarray*}
For any $x \in \mathcal{V}$ and as $\sharp{\widetilde{K}}< \infty$, from Lemma (\ref{lem}), we get
\begin{equation}\label{d}
f^{2}(x) \leq C_{1}\left(\norm{f}_{l^{2}({\widetilde{K}})}^{2}+ \norm{\mathrm{d}f }_{l^{2}(\mathcal{E})}^{2}\right).
\end{equation}
But, $\norm{f}_{l^{2}({\widetilde{K}})}=0$ and $\norm{\mathrm{d}f }_{l^{2}(\mathcal{E})}=0$. So it follows immediately that $f=0$ on $\mathcal{V}$.\\

It remains to show that if $\norm{ \varphi}_{l^{2}(\mathcal{E}_{\widetilde{K}})}=0$ and $\norm{\delta \varphi }_{l^{2}(\mathcal{V})}=0 $ then $\varphi=0$. We suppose that $\varphi \neq 0$.

$\varphi$ is a finite support function in $\mathcal{E}\setminus\mathcal{E}_{\widetilde{K}}$ and therefore, by Lemma (\ref{lem}) where $U$ equals to the support of $\varphi$, there exists a positive constant $C$ such that
$$C \norm{ \varphi}_{l^{2}(\mathcal{E}_{U})}\leq \norm{ \varphi}_{l^{2}(\mathcal{E}_{\widetilde{K}})}+\norm{\delta \varphi }_{l^{2}(\mathcal{V})}.$$
But, $\norm{ \varphi}_{l^{2}(\mathcal{E}_{\widetilde{K}})}=\norm{\delta \varphi }_{l^{2}(\mathcal{V})}=0$, since we get $\varphi = 0$ on $\mathcal{E}_{U}$, which is impossible.\\

\underline{\textbf{Aim ii)}} Show that the space $W$ is independent of the choice of $G_{\widetilde{K}}$.\\

Let $G_{\widetilde{K}_{1}}$ be another neighborhood of $G_{K}$ such that $K\subset \widetilde{K}_{0} \subset \widetilde{K}_{1}$.\\

So, we have
$$N_{\widetilde{K}_{0}}(f, \varphi)\leq N_{\widetilde{K}_{1}}(f, \varphi).$$

Moreover, to show the existence of a constant $C>0$ such that $N_{\widetilde{K}_{1}}(f, \varphi)\leq C N_{\widetilde{K}_{0}}(f, \varphi)$, it suffices to show the existence of a constant $C>0$ such that $\norm{(f, \varphi)}_{l^{2}(\widetilde{K}_{1}\setminus\widetilde{K}\hat{_{0}})}^{2}\leq C N_{\widetilde{K}_{0}}^{2}(f, \varphi)$. Indeed, we have
\begin{eqnarray*}
N_{\widetilde{K}_{1}}^{2}(f, \varphi)&=&\norm{(f, \varphi)}_{l^{2}(\widetilde{K}_{1})}^{2}+\norm{D(f, \varphi)}_{l^{2}(G)}^{2}\\
&=& \norm{(f, \varphi)}_{l^{2}(\widetilde{K}_{1}\setminus\widetilde{K}_{0})}^{2}+\norm{(f, \varphi)}_{l^{2}(\widetilde{K}_{0})}^{2}+\norm{D(f, \varphi)}_{l^{2}(G)}^{2}\\
&=& \norm{(f, \varphi)}_{l^{2}(\widetilde{K}_{1}\setminus\widetilde{K}_{0})}^{2}+N_{\widetilde{K}_{0}}^{2}(f, \varphi).
\end{eqnarray*}

Using lemma (\ref{lem}) and as we have $\sharp ({\widetilde{K}_{1}}\setminus{\widetilde{K}_{0}})<\infty$, we get

\begin{equation*}\norm{f}_{l^{2}({\widetilde{K}_{1}}\setminus{\widetilde{K}_{0}})}^{2} \leq C \left(\norm{f}_{l^{2}({\widetilde{K}_{0}})}^{2}+ \norm{\mathrm{d}f }_{l^{2}(\mathcal{E})}^{2}\right),\end{equation*}

where $C=C({\widetilde{K}_{1}}\setminus{\widetilde{K}_{0}}, {\widetilde{K}_{0}})$.\\

And

\begin{equation*}\norm{\varphi}_{l^{2}(\mathcal{E}_{\widetilde{K}_{1}}\setminus\mathcal{E}_{\widetilde{K}_{0}})}^{2} \leq C \left(\norm{\varphi}_{l^{2}(\mathcal{E}_{\widetilde{K}_{0}})}^{2}+ \norm{\delta \varphi }_{l^{2}(\mathcal{V})}^{2}\right).\end{equation*}

where $C=C({\widetilde{K}_{1}}\setminus{\widetilde{K}_{0}}, {\widetilde{K}_{0}})$.\\

So, we obtain

$$\norm{(f, \varphi)}_{l^{2}(G_{\widetilde{K}_{1}}\setminus G_{\widetilde{K}_{0}})}^{2}\leq C N_{\widetilde{K}_{0}}^{2}(f, \varphi).$$

Thus, we have shown that the construction of a norm on $W$ is independent of the choice of the neighborhood associated to the subgraph $G_K$. We set:

$$\norm{(f, \varphi)}_{W}:=\left(\norm{(f, \varphi)}_{l^{2}(G_{\widetilde{K}_{0}})}^{2}+\norm{D(f, \varphi)}_{l^{2}(G)}^{2}\right)^{\frac{1}{2}},$$
for $(f, \varphi)\in \mathcal{C}_{0}(\mathcal{V}) \oplus \mathcal{C}^{a}_{0}(\mathcal{E}) $.\\

\underline{\textbf{Aim iii)}}: By Lemma (\ref{lem}), we have the injection of $\mathcal{C}_{0}(\mathcal{V}) \oplus \mathcal{C}^{a}_{0}(\mathcal{E})$ to
$\mathcal{C}(\mathcal{V}) \oplus \mathcal{C}^{a}(\mathcal{E})$ extends by continuity to $W$.

\underline{\textbf{Aim iv)}}: we have
$$\norm{D(f, \varphi)}_{l^{2}(G)}^{2}\leq \norm{(f, \varphi)}_{l^{2}(G_{\widetilde{K}})}^{2} + \norm{D(f, \varphi)}_{l^{2}(G)}^{2}=\norm{(f, \varphi)}_{W}^{2}.$$

Consequently, $D:W \longrightarrow l^{2}(G)$ is a bounded operator.

\end{Demo}

\section{Semi-Fredholmness of the discrete Gauss-Bonnet operator}
\begin{defi}
An operator is semi-Fredholm if its range is closed and its kernel is finite dimensional .
\end{defi}
Now we come to our main result:\\

\textbf{Theorem.}
Let $W$ be a Hilbert space satisfying:
\begin{enumerate}
  \item$\mathcal{C}_{0}(\mathcal{V}) \oplus \mathcal{C}^{a}_{0}(\mathcal{E})$ is dense in $W$.

\item The injection of $\mathcal{C}_{0}(\mathcal{V}) \oplus \mathcal{C}^{a}_{0}(\mathcal{E})$ to $\mathcal{C}(\mathcal{V}) \oplus \mathcal{C}^{a}(\mathcal{E})$  extends by continuity to $W$.

  \item $D:W \longrightarrow l^{2}(G)$ is a bounded operator.
  \end{enumerate}

  Then, the following conditions are equivalent:\\
\;\;\;  i) $D:W \longrightarrow l^{2}(G)$ \;is\; semi-Fredholm.\\
 \;\; ii) There exists a finite subgraph $G_{K}$ of $G$ and a positive constant $C=C_{K}$ such that
\begin{equation}
C \norm{(f, \varphi)}_{W} \leq \norm{D(f, \varphi)}_{l^{2}(G)},\;\forall (f, \varphi) \in \mathcal{C}_{0}(\mathcal{V}\setminus K)\times \mathcal{C}^{a}_{0}(\mathcal{E}\setminus \mathcal{E}_{K}). \end{equation}

\begin{Demo}

We take the same arguments used by Carron \cite{CARR}. We start by showing the direct implication, we assume that the conclusion is false. Then, we can find an increasing sequence of finite subgraph $\{G_{K_{n}}\}_{n}$ such that $G=\bigcup_{n} G_{K_{n}}$ and a sequence $\{\sigma_{n}\}_{n}$ with finite support in $\mathcal{V}\setminus K_{n}$ satisfying the following conditions, for all $n\geq 1$

$$\left\{
  \begin{array}{ll}\sigma_{n}=(f_{n}, \varphi_{n})\in \mathcal{C}_{0}(\mathcal{V}\setminus K_{n})\times \mathcal{C}^{a}_{0}(\mathcal{E}\setminus\mathcal{E}_{K_{n}}),\\
\\
   \norm{\sigma_{n}}_{W}=1,\\
\\
\norm{D\sigma_{n}}_{l^{2}(G)}\leq\frac{1}{n} .
\end{array}
\right.
$$

On the other hand, it was assumed that $D:W\longrightarrow l^{2}(G)$ is semi-Fredholm. Therefore, by \cite{shubin} there exists a bounded operator $P:l^{2}(G)\longrightarrow W$ such that

\begin{equation}\label{P inv}
P\circ D =\mathrm{Id}_{W}-H,\end{equation} where $H$ is the orthogonal projection onto the kernel of $D$, it is an operator with finite rank. 

Then, we obtain
\begin{eqnarray*}
\norm{\sigma_{n}}_{W}&\leq&\norm{(P\circ D)\sigma_{n}}_{W}+\norm{H\sigma_{n}}_{W}\\
&\leq&\norm{P}\norm{ D\sigma_{n}}_{l^{2}(G)}+\norm{H\sigma_{n}}_{W}\\
&\leq&\left(\frac{\norm{P}}{n}+\norm{H\sigma_{n}}_{W}\right).
\end{eqnarray*}
If $$ \lim\limits_{n \to \infty}\norm{H\sigma_{n}}_{W} = 0 \Longrightarrow \lim\limits_{n \to \infty}\norm{\sigma_{n}}_{W} = 0, $$
which contradicts the assumption $\norm{\sigma_{n}}_{W}=1$.\\

So, our aim is to prove that $ \{H\sigma_{n}\}_{n}$ converges to $0$ in $W$. Indeed, we set

\begin{equation}\label{sig}\sigma_{n}=\sigma_{n}^{1}+\sigma_{n}^{2}\end{equation}

with $\sigma_{n}^{1}(=H\sigma_{n})\in \ker D$ and $\sigma_{n}^{2}\in (\ker D)^{\bot}$.\\

Such as
$$\left\{
  \begin{array}{ll}(P\circ D)\sigma_{n}=\sigma_{n}^{2},\\
\\
\norm{P\circ D\sigma_{n}}_{W}\leq \norm{P}\norm{D\sigma_{n}}_{l^{2}(G)}\longrightarrow_{n\rightarrow \infty} 0.
\end{array}
\right.
$$
Then, for the norm of $W$

\begin{equation}\label{limit}\lim_{n\rightarrow \infty} \sigma_{n}^{2}= 0.\end{equation}

Moreover, $\{\sigma_{n}^{1}\}_n$ is a bounded sequence of $\ker D$ which is of finite dimension. So we can extract a subsequence converging  to $\sigma$ in $W$, which we denote $\{\sigma_{\varphi(n)}^{1}\}_n$.\\

Using (\ref{sig}) and (\ref{limit}), $\{\sigma_{\varphi(n)}\}_n$ converges in $W$ to $\sigma$ (as a sum of two converging sequences) and as a consequence $\norm{\sigma}_{W}=1$.\\  

Let us prove that $\sigma=0$ where $\sigma=\lim \sigma_{\varphi(n)}=\lim\sigma_{\varphi(n)}^{1}$.\\

We suppose that $\sigma\neq 0$. As $W$ is injected continuously in $\mathcal{C}(\mathcal{V}) \oplus \mathcal{C}^{a}(\mathcal{E})$ , there exists $x\in V$ such that $\{\sigma_{\varphi(n)}(x)\}_{n}$ converges to $\sigma(x)\neq 0$. But, by construction the sequence $\{\sigma_{\varphi(n)}\}_n$ converges ponctually to $0$ ( the sequence $\{\sigma_{\varphi(n)}\}_n$ has a finite support outside of $G_{K_{n}}$). Hence, we conclude that $\sigma(x)=0$ which is absurd.\\





It remains to prove $ii) \Rightarrow i)$.\\

\textbf{First step:} We construct a bounded operator $Q:l^{2}(G)\longrightarrow W$ such that $Q\circ D-\mathrm{Id}_{W}$ is a compact operator, this will show that $D:W\longrightarrow l^{2}(G)$ has a finite kernel and a closed range \cite{shubin}.\\

Let $D_{1}$ be the restriction of $D$ on $G\setminus G_{K}$, so $D_{1}:W(G\setminus G_{K})\longrightarrow l^{2}(G)$ is bounded, where $W(G\setminus G_{K})=\{\sigma=(f, \varphi) \in W;\; \sigma=0 \hbox{ on }G_{K}\}$. Moreover, by assumption we have
$$C \norm{(f, \varphi)}_{W} \leq \norm{D(f, \varphi)}_{l^{2}(G)},\;\forall (f, \varphi) \in \mathcal{C}_{0}(\mathcal{V}\setminus K)\times \mathcal{C}^{a}_{0}(\mathcal{E}\setminus\mathcal{E}_{K}). $$
Then, $D_{1}$ is injective with closed range, which allows the existence of a left inverse $P_{1}$ such that
$$P_{1}\circ D_{1}=\mathrm{Id}.$$
On the other hand, we denote $$D_{2}:l^{2}(\widetilde{K}_{1}) \longrightarrow l^{2}(G) $$
where $\widetilde{K}_{1}$ is a neighborhood (see Definition $(\ref{neighb})$) of $\widetilde{K}_{0}$, such that $\widetilde{K}_{0}$ is the smallest neighborhood of $K$.\\

Since $l^{2}(\widetilde{K}_{1})$ is a vector space of finite dimension, then $D_{2}$ is continuous with closed range. We denote $P_{2}$ "the parametrix" which is a continuous operator satisfying
$$P_{2}\circ D_{2}=\mathrm{Id}-H_{2},$$
where $H_{2}$ is the orthogonal projection onto the kernel of $D_{2}$.\\

We consider now the indicator function $\chi$ as in $(\ref{cut})$ by replacing $K$ by $\widetilde{K}_{0}$, which gives $\mathrm{d}\chi$, $\overline{\chi}$, $1-{\chi}$ and $1-\overline{\chi}$ where
$$(1-\chi)(x)=\left\{
  \begin{array}{ll}
  1  \hbox{ if }  x\in \widetilde{K}_{0}\\
\\
   0  \hbox{ otherwise}.
   \end{array}
\right.
 \hbox{ and }(1-\overline{\chi})(e)=\left\{
  \begin{array}{ll}
  1 \;\hbox{if}\;e\in \mathcal{ E}_{\widetilde{K}_{0}},\\
\\
   \frac{1}{2}\; \hbox{if} \;e\in \partial \mathcal{E}_{\widetilde{K}_{0}},\\
   \\
   0\; \hbox{otherwise}.
   \end{array}\right.$$

Furthermore, we define the operator $\chi .$ depending on the domain by:\\

If  $\chi .:C_{0}(\mathcal{V})\longrightarrow C_{0}(\mathcal{V})$ so we have $\chi .f= \chi f$, for all $f\in C_{0}(\mathcal{V})$.\\

If  $\chi .:C^{a}_{0}(\mathcal{E})\longrightarrow C^{a}_{0}(\mathcal{E})$  we get $\chi .\varphi= \overline{\chi} \varphi$, for all $\varphi\in C^{a}_{0}(\mathcal{E})$.\\

If  $\chi .:C_{0}(\mathcal{V})\oplus C^{a}_{0}(\mathcal{E}) \longrightarrow C_{0}(\mathcal{V})\oplus C^{a}_{0}(\mathcal{E})$ hence we obtain $\chi .(f, \varphi)= (\chi f,\overline{\chi} \varphi)$,\\

for all $(f, \varphi)\in C_{0}(\mathcal{V})\oplus C^{a}_{0}(\mathcal{E})$.\\

We set
$$Q\sigma:=P_{2}(1-{\chi})\sigma+P_{1}{\chi}\sigma,$$
where $\sigma=(f, \varphi)$.\\

\textbf{Second step:} Let us check that the operator $Q\circ D-\mathrm{ Id}$ is compact.\\
We denote the following bracket for any two operators $A$ and $B$:
$$[A, B]=A B-B A.$$
Then, we obain
\begin{eqnarray*}
Q\circ D&=&P_{2}(1-{\chi})D+P_{1}{\chi}D\\
&=&P_{2}D(1-{\chi})+  P_{2}[1-{\chi}, D]  + P_{1}D{\chi}+  P_{1}[{\chi}, D]\\
&=&P_{2}D_{2}(1-{\chi})+  P_{2}[1-{\chi}, D]  + P_{1}D_{1}{\chi}+  P_{1}[{\chi}, D]\\
&=&(\mathrm{Id}-H_{2})(1-{\chi})+  P_{2}[1-{\chi}, D]  + \mathrm{Id}({\chi})+  P_{1}[{\chi}, D]\\
&=& \mathrm{Id}-H_{2}(1-{\chi})+  P_{2}[1-{\chi}, D]   +  P_{1}[{\chi}, D].
\end{eqnarray*}
We just calculate $P_{2}[{\chi}, D]$. We have
$$[{\chi}, D]=[{\chi}, \mathrm{d}]+[{\chi}, \delta]. $$
For the first bracket, we obtain
\begin{eqnarray*}
[{\chi}, \mathrm{d}]f(e)&=&\overline{\chi}(e) \mathrm{d}(f)(e) - \mathrm{d}({\chi}f)(e)\\
&=&\frac{1}{2}\left[{\chi}(e^{+})+{\chi}(e^{-}) \right] \mathrm{d}(f)(e) -{\chi}(e^{+}) \mathrm{d}(f)(e)-f(e^{-}) \mathrm{d}{\chi}(e)\\
&=&-\frac{1}{2} \mathrm{d}{\chi}(e) \mathrm{d}(f)(e)-f(e^{-}) \mathrm{d}{\chi}(e).
\end{eqnarray*}
And for the second one, we get
\begin{eqnarray*}
[{\chi}, \delta]\varphi(x)&=&{\chi}(x) \delta(\varphi)(x) -\delta(\overline{\chi}\varphi)(x)\\
&=&{\chi}(x) \delta(\varphi)(x)-{\chi}(x) \delta(\varphi)(x)+\frac{1}{2} \sum_{e, e^{+}= x}\mathrm{d}({\chi})(e)\varphi(e)\\
&=&\frac{1}{2} \sum_{e, e^{+}= x}\mathrm{d}({\chi})(e)\varphi(e).
\end{eqnarray*}

But, the support of $\mathrm{d}({\chi})$ is included in $\partial \mathcal{E}_{\widetilde{K}_{0}}\subset \widetilde{K}_{1}$ which is finite. Then, $[{\chi}, D] $ has a finite range so it is a compact operator.

Finally, $Q\circ D=\mathrm{Id}+H$ where $H$ is a compact operator .\\

\end{Demo}
\begin{rema}In the Theorem , we obtain $D$ Fredholm if it is an essential-selfadjoint operator \cite{CARR}. \end{rema}
\begin{rema}
There is a second method inspired from \cite{AN} to show ii) $\Rightarrow$ i) of the Theorem . This can be demonstrated with the aid of the following \textbf{claim}:
"If $\sigma_{n}=(f_{n}, \varphi_{n})\in \mathcal{C}_{0}(\mathcal{V})\times \mathcal{C}^{a}_{0}(\mathcal{E})$ is $W$-bounded and $(D\sigma_{n})_{n}$ is convergent in $l^{2}(G)$, then $(\sigma_{n})_{n}$ has a $W$-convergent subsequence".\\
\end{rema}
We have the following result:
\begin{Prop} Let $W$ be a Hilbert space satisfying:

\begin{enumerate}
  \item$\mathcal{C}_{0}(\mathcal{V}) \oplus \mathcal{C}^{a}_{0}(\mathcal{E})$ is dense in $W$.

\item The injection of $\mathcal{C}_{0}(\mathcal{V}) \oplus \mathcal{C}^{a}_{0}(\mathcal{E})$ to $\mathcal{C}(\mathcal{V}) \oplus \mathcal{C}^{a}(\mathcal{E})$  extends by continuity to $W$.

  \item $D:W \longrightarrow l^{2}(G)$ is a bounded operator.
  \end{enumerate}
Then if there exists a finite subgraph $G_{K}$ of $G$ and a positive constant $C=C_{K}$ such that

\begin{equation}\label{prop1}
C \norm{(f, \varphi)}_{W} \leq \norm{D(f, \varphi)}_{l^{2}(G)},\;\forall (f, \varphi) \in \mathcal{C}_{0}(\mathcal{V}\setminus K)\times \mathcal{C}^{a}_{0}(\mathcal{E}\setminus \mathcal{E}_{K}), \end{equation}

so necessarily, the operator $D:W \longrightarrow l^{2}(G)$ is semi-Fredholm.
\end{Prop}

\begin{Demo}

We start by proving the following \textbf{claim}: if $\sigma_{n}=(f_{n}, \varphi_{n})\in \mathcal{C}_{0}(\mathcal{V})\times \mathcal{C}^{a}_{0}(\mathcal{E})$ is $W$-bounded and $(D\sigma_{n})_{n}$ is convergent in $l^{2}(G)$, then $(\sigma_{n})_{n}$ has a $W$-convergent subsequence.\\

Let $G_{\widetilde{K}}$ be a neighborhood of the subgraph $G_{K}$ (see Definition \ref{neighb}), then $({\sigma_{n}}\upharpoonright_{\widetilde{K}})_{n}$ is a bounded sequence in a vector space with finite dimension. Hence, it admits a convergent subsequence.\\

In $G\setminus G_{\widetilde{K}}$, we consider the indicator function $\chi$ as in $(\ref{cut})$ by replacing $K$ by $\widetilde{K}$. Then, we obtain a function $\chi\sigma_{n}$ with finite support in $G\setminus G_{K}$ and we can apply the inequality (\ref{prop1}) to $\chi\sigma_{n}$, in particular to $(\chi f_{n}, 0)$ and $(0, \overline{\chi}\varphi_{n})$. First, we obtain
$$\norm{\chi f_{n}}_{W}\leq C\norm{\mathrm{d}(\chi f_{n})}_{l^{2}(\mathcal{E})}.$$
But, from the equality (\ref{eqlem3}) of Lemma (\ref{lem3}), we get
$$\mathrm{d}(\chi f_{n})(e)=\chi(e^{+})\mathrm{d}(f_{n})(e)+f_{n}(e^{-})\mathrm{d}(\chi)(e).$$
We have $(d(f_{n}))_{n}$ is a convergent sequence and  $supp(d\chi)\subset \mathcal{E}_{\widetilde{K}}$ is finite, thus,
$f_{n}(x)\upharpoonright_{{\widetilde{K}}}$ admits a convergent subsequence. \\
Then we may conclude that $\chi f_{n}$ admits a $W$-convergent subsequence, i.e, $(f_{n}\upharpoonright_{V\setminus \widetilde{K}})_{n}$ admits a $W$-convergent subsequence.\\

Second, we have
$$\norm{\overline{\chi}\varphi_{n}}_{W}\leq C\norm{\delta(\overline{\chi}\varphi_{n})}_{l^{2}(\mathcal{V})}.$$
Since the equality (\ref{eqqlem3}) of Lemma (\ref{lem3}) gives
$$\delta(\overline{\chi}\varphi_{n})(x)=\chi(x)\delta(\varphi_{n})(x)-\frac{1}{2c(x)}\sum_{e, e^{+}= x}r(e)\mathrm{d}(\chi)(e) \varphi_{n}(e), \;\forall x\in \mathcal{V}.$$
Furthermore by assumptions the sequence $(\delta(\varphi_{n}))_{n}$ is convergent and $\mathrm{supp}(\mathrm{d}\chi)\subset \mathcal{E}_{\widetilde{K}}$ is finite, hence, $(\varphi_{n}\upharpoonright_{\mathcal{E}_{\widetilde{K}}})$ admits a convergent subsequence.
As a result, we deduce that the sequence $(\overline{\chi}\varphi_{n})_{n}$ admits a $W$-convergent subsequence. So, the sequence $(\varphi_{n}\upharpoonright_{\mathcal{E}\setminus \mathcal{E}_{\widetilde{K}}})_{n}$ admits a $W$-convergent subsequence.\\

Now we can show that our operator $D$ is semi-Fredholm.
\begin{enumerate}
  \item We start by proving that $\ker D$ is finite dimensional, which is equivalent to show that
   $\{ \sigma \in \ker D;\; \norm{\sigma}_{W}=1 \}$ is compact.\\

   Let $(\sigma_{n})_{n}\subset \ker D$ be such that $\norm{\sigma_{n}}_{W}=1$ and $D\sigma_{n}=0$. Then, by the claim,
$(\sigma_{n})_{n}$ admits a convergent subsequence. So the result occurs.

  \item Let us show that ${\Ima}D$ is closed.\\

 Let $(y_{n})_{n}$ be a sequence of ${\Ima}D$ such that $(y_{n})_{n}$ converges to $y$ in
  $l^{2}(G)$. Is that $y$ in ${\Ima}D$? \\

  Since $(y_{n})_{n} \subset {\Ima}D $, then there exist
 $(\sigma_{n})_{n}\subset \ker D^{\bot}$ and $\sigma_{n}\neq0 \;\forall n$,
  such that $y_{n}=D\sigma_{n}$.
 $(\sigma_{n})_{n}$ must be bounded. If not, by extraction we can construct $s_{n}=\frac{\sigma_{n}}{\norm{\sigma_{n}}_{W}}$, such that
  $$\left\{
  \begin{array}{ll}
 ( s_{n})_{n} \subset \ker D^{\bot}\\
\\
   \norm{s_{n}}_{W}=1 \\
   \\
   Ds_{n}\rightarrow 0.
   \end{array}
\right.$$
Using the claim, we can conclude that $( s_{n})_{n}$ admits a convergent subsequence with limit denoted $s$ such that
  $$\left\{
  \begin{array}{ll}
  s \in \ker D^{\bot}\\
\\
   \norm{s}_{W}=1 \\
   \\
   Ds= 0.
   \end{array}
\right.$$
Then, $s \in \ker D \cap \ker D^{\bot}=\{0\}$. So $s=0$, which is absurd.\\

Hence the sequence $( \sigma_{n})_{n}$ is bounded and since $( D\sigma_{n})_{n}$ converges to $y$, using the claim, the sequence $( \sigma_{n})_{n}$ admits a convergent subsequence and let $\sigma$ be this limit. But, the operator $D$ is bounded.
Then, $D\sigma_{n}$ converges to $D\sigma$ and by uniqueness of the limit $y=D{\sigma}$.

\end{enumerate}

\end{Demo}

\begin{Corol}
$D$ is non-parabolic at infinity if and only if there exists a finite subgraph $G_{K}$ of $G$ such that if we complete $\mathcal{C}_{0}(\mathcal{V})\times C^{a}_{0}(\mathcal{E})$ by the norm
$$\norm{(f, \varphi)}_{W}=\left(\norm{(f, \varphi)}_{l^{2}(\widetilde{K})}^{2}+\norm{D(f, \varphi)}_{l^{2}(G)}^{2}\right)^{\frac{1}{2}},$$
in order to obtain $W$ satisfying \begin{enumerate}
  \item$\mathcal{C}_{0}(\mathcal{V}) \oplus \mathcal{C}^{a}_{0}(\mathcal{E})$ is dense in $W$.

\item The injection of $\mathcal{C}_{0}(\mathcal{V}) \oplus \mathcal{C}^{a}_{0}(\mathcal{E})$ to $\mathcal{C}(\mathcal{V}) \oplus \mathcal{C}^{a}(\mathcal{E})$  extends by continuity to $W$.

  \item $D:W \longrightarrow l^{2}(G)$ is semi-Fredholm.
  \end{enumerate}
\end{Corol}

\section{Examples}
\subsection{A star-like graph}
\begin{defi}
The disjoint union of two graphs $G_\alpha=(\mathcal{V}_\alpha, \mathcal{E}_\alpha)$ and $G_\beta=(\mathcal{V}_\beta, \mathcal{E}_\beta)$ is the disjoint union of their vertex and edge with no edge joining $\mathcal{V}_\alpha$ and $\mathcal{V}_\beta$.
\end{defi}
According to \cite{ytt1}, we have the following definition:

\begin{defi}
An infinite graph $G=(\mathcal{V}, \mathcal{E})$ is called \textbf{star-like}, if there exists a finite subgraph $G_K$ of $G$ so that $G\setminus G_K$ is the union of a finite number of disjoint copies $G_\alpha$ of the graph $\N$.
\end{defi}

\begin{figure}[ht]\label{}
\begin{center}
\includegraphics*[scale=0.7]{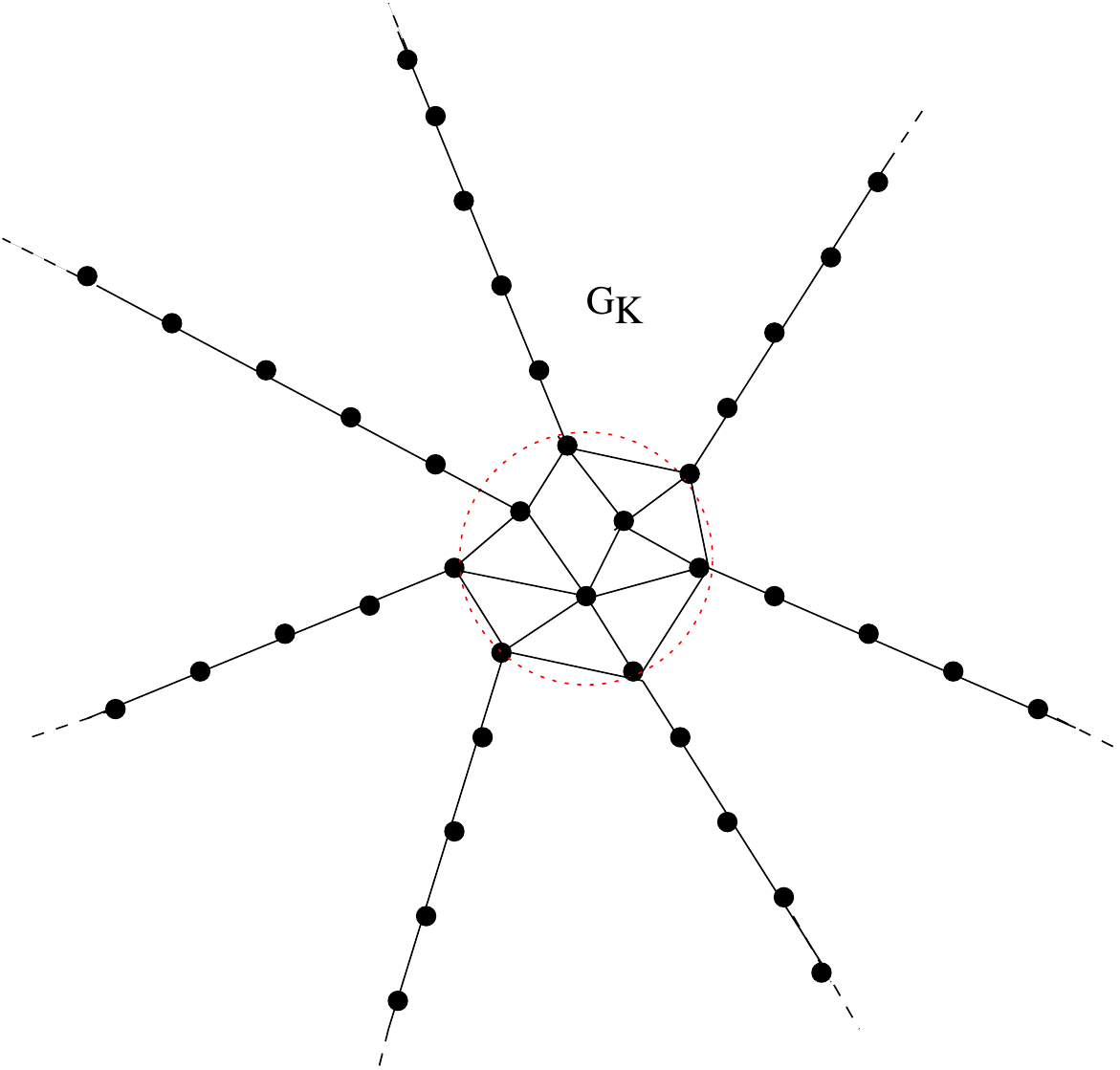}
\vspace{-0.5cm}
\caption{}
\end{center}
\end{figure}

\begin{Prop}
 In the case where $c=r=1$, $D$ is non-parabolic at infinity in the star-like graph.
\end{Prop}
\begin{Demo}

By the definition of the star-like graph, there exists a finite subgraph $G_{K}$ of $G$ so that $G\setminus G_K=\bigsqcup_{\alpha \in J} G_\alpha$. Let $U$ be a finite subset of $G\setminus G_K$ then, there exists $\alpha \in J$ such that $U\subset G_\alpha$. We look for a positive constant $C=C(U)$ such that
\begin{equation}\label{eq exemp}C \norm{(f, \varphi)}_{l^{2}(U)} \leq \norm{D(f, \varphi)}_{l^{2}(G)},\;\forall (f, \varphi)\in \mathcal{C}_{0}(\mathcal{V}\setminus K) \oplus \mathcal{C}^{a}_{0}(\mathcal{E}\setminus \mathcal{E}_{K}).\end{equation}\\
Let $f \in \mathcal{C}_{0}(\mathcal{V}\setminus K)$ such that $U$ is included in the support of $f$.\\

For $U=\{a\}$, we have
$$\norm{f}_{l^{2}(U)}^{2}=f^{2}(a).$$
For $o \in K$ and as $G$ is connected we can find a path $ \gamma_{oa}$ joining $o$ to $a$. Suppose that this path is of length $n$ such that $x_{0}=a$ and $x_{n}=o$, using the Jensen's inequality and $f(x_{n})=0$, we obtain
\begin{eqnarray*}
f^{2}(a)&=&\left( f(x_{0})-f(x_{1})+ f(x_{1})-f(x_{2})+ f(x_{2})-...-f(x_{n-1})+ f(x_{n-1})-f(x_{n})+ f(x_{n}) \right)^{2}\\
&\leq& n\left(\left(f(a)-f(x_{1})\right)^{2}+\left(f(x_{1})-f(x_{2})\right)^{2}+...+\left(f(x_{n-1})-f(x_{n})\right)^{2}\right),
\end{eqnarray*}
which implies
\begin{equation}\label{exemp}f^{2}(a)\leq n \norm{\mathrm{d}f}_{l^{2}(V)}^{2}.\end{equation}
\begin{rema} $n$ depends only on $U$ and $K$.\end{rema}
Similarly, for $\varphi \in C^{a}_{0}(\mathcal{E} \setminus \mathcal{E}_{K})$, we obtain
$$\norm{\varphi}_{l^{2}(\mathcal{E}_{U})}^{2}\leq C_{U} \norm{\delta \varphi}_{l^{2}(\mathcal{V})}^{2}.$$
Moreover, for $U=\{a_{1},...,a_{n}\}$, we prove the inequality (\ref{eq exemp}).\\

By the inequality (\ref{exemp}), for all $i \in \{1,...,n\}$, we get
$$f^{2}(a_{i})\leq {n_{i}}\norm{\mathrm{d}f}_{l^{2}(\mathcal{V})}^{2}$$
where $n_{i}$ is the number of edge of the shortest path between $a_{i}$ and any vertex of $K$.\\

For thus, we have
$$\sum_{i=1}^{n} f^{2}(a_{i})\leq \sum_{i=1}^{n}{n_{i}} \norm{\mathrm{d}f}_{l^{2}(\mathcal{V})}^{2}.$$
Hence
$$\norm{f}_{l^{2}(U)}^{2}\leq C_{U} \norm{\mathrm{d}f}_{l^{2}(\mathcal{V})}^{2}.$$
And similarly, we show that
$$\norm{\varphi}_{l^{2}(\mathcal{E}_{U})}^{2}\leq C_{U} \norm{\delta \varphi}_{l^{2}(\mathcal{V})}^{2}.$$

\end{Demo}

\subsection{The triadic tree}
\begin{defi}
A tree is a connected graph containing no cycles. The \textbf{triadic tree} is the tree such that all the vertices have degree 3.
\end{defi}

\begin{figure}[ht]\label{}
\begin{center}
\includegraphics*[scale=0.8]{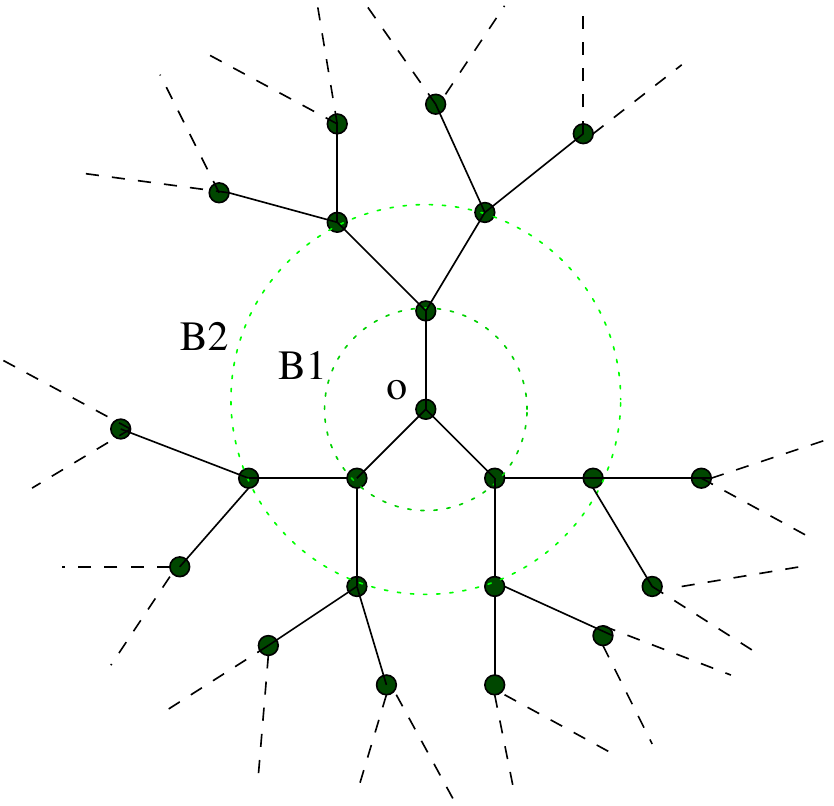}
\vspace{-0.25cm}
\caption{}
\end{center}
\end{figure}

\begin{Prop}
In the triadic graph the condition of "non-parabolicity at infinity" is not verified.
\end{Prop}

\begin{Demo}

We fix a vertex $o$, see the figure 2, we can find an increasing sequence of finite subgraph $\{G_{n}\}_{n}$ such that  $G_{n}=\{x\in \mathcal{V};\;d(o, x)\leq n\}$ and $G=\bigcup_{n} G_{n}$.
The contradiction of non-parabolicity at infinity property could be: for all $n$ there exists $U$ outside of $G_{n}$ and a 1-form $\varphi_{n}$ with finite support outside of $G_{n}$
such that $\delta \varphi_{n}=0$ and $\norm{\varphi_{n}}_{l^{2}(U)} \neq 0$. Such $\varphi_{n}$ exist. Indeed one can construct a skewsymmetric function $\varphi_{n}$ supported on the outward tree of every vertex $x_{n} \in G_{n}$  with $\delta \varphi_{n}=0$ in the following way: let $e_{0}$ and $b_{0}$ be the two outward edges of $x_{n}$ (the third one rely $x_{n}$ to $x_{n-1}$) and denote $e_{m}^{k}$, $m\geqslant 1$, $1\leqslant k \leqslant 2 ^{m}$, resp. $b_{m}^{k}$, $m\geqslant 1$, $1\leqslant k \leqslant 2 ^{m}$, the outward edges emanating from $e_{0}$, resp. $b_{0}$, of generation $m$. We define $\varphi_{n}$ to be $0$ excepted on these edges where $\varphi_{n}(e_{m}^{k})=\frac{1}{2 ^{m}}$ and $\varphi_{n}(b_{m}^{k})=-\frac{1}{2 ^{m}}$ (the edge are oriented outward). So, we deduce that $\delta$ does not satisfy the property of non-parabolicity at infinity.

\end{Demo}

\begin{rema}
We can generalize this example for the tree with degree $d\geq 3$, we can use the same argument with $ \varphi_{n} =\pm(\frac{1}{d-1})^m$.
\end{rema}

\begin{rema}

a) The importance of non-parabolicity at infinity appears with the operator $\delta$. In fact, this property for the operator $\mathrm{d}$ is always true on any connected graph.\\
b) In probability \cite{gry} and potential theory \cite{wo} there exists an interesting notion of non-parabolic for the graph which is equivalent (\cite{A} Theorem 2.1) to the following statement: there exists $\;x \in \mathcal{V}$ \;and\; $C > 0$ such that $$f^{2}(x)\leq C \norm{\mathrm{d}f}^{2}_{l^{2}(\mathcal{E})},\; \forall f\in \mathcal{C}_{0}(G).$$ This notion is different from the non-parabolicity at infinity. Indeed, the graph $\mathbb{Z}$ and $\mathbb{Z}^{2}$ are parabolic, but $\mathbb{Z}^{n}$, $n\geq 3$ is non-parabolic. On the other side, we have $\delta$ is non-parabolic at infinity in $\mathbb{Z}$ but in $\mathbb{Z}^{n}$, $n\geq 2$, $\delta$ does not verify this property (since it has cycles supported outside any finite subgraph).

\end{rema}

\textbf{Acknowledgements}
I would like to express my sincere gratitude to my advisors Professors Nabila Torki-Hamza and Colette Anné for the continuous support of my Ph.D study, for their patience, motivation, enthusiasm, and immense knowledge. Their guidance helped me in all the time of research and without their wisdom this paper would not have been possible. I would like to thank the Laboratory of Mathematics Jean Leray of Nantes (LMJL) and the research unity (UR/13ES47) of Faculty of Sciences of Bizerte (University of Carthage) for  its  financial and its continuous support. Also, this work was financially supported by the "PHC Utique" program of the French Ministry of Foreign Affairs and Ministry of higher education and research and the Tunisian Ministry of higher education and scientific research in the CMCU project number 13G1501 " Graphes, Géométrie et théorie Spectrale".\\
Finally, I take this chance to thank Jun Masamune for the fruitful discussions during his visit to LMJL Nantes. I would like to thank also the anonymous referee for the careful reading of my paper and the valuable comments and suggestions.

\vskip -2.5 cm

\end{document}